\documentclass[a4paper]{amsart}
\usepackage{paralist}
\usepackage{hhline}

\usepackage{longtable}
\usepackage{amssymb, amsthm}

\usepackage[dvips, dvipsnames, usenames]{color}

\newcommand{\fd}{finite-dimensional}

\makeatletter

\newcommand{\Aut}{\operatorname{Aut}}
\newcommand{\End}{\operatorname{End}}
\newcommand{\Irr}{\operatorname{Irr}}
\newcommand\id{\operatorname{id}}

\newcommand\sgn{\operatorname{sgn}}

\newcommand{\atlas}{\textsf{ATLAS}}

\newcommand{\ku}{\Bbbk}

\newcommand{\GL}{\mathbf{GL}}

\newcommand{\ydg}{{}^{\ku G}_{\ku G}\mathcal{YD}}

\newcommand\ac{\mathbb A_4}

\newcommand{\Fc}{{\mathcal F}}
\newcommand{\oc}{{\mathcal O}}
\newcommand{\Oc}{{\mathcal O}}
\newcommand\s{\mathbb S}

\newcommand\sn{\mathbb S_n}

\newcommand\bn{\mathbb B_n}

\newcommand{\N}{{\mathbb N}}

\providecommand{\\}{\\}

\numberwithin{equation}{section} 
  \theoremstyle{plain}
  \newtheorem{thm}{Theorem}[section]
  \theoremstyle{plain}
  \newtheorem{lem}[thm]{Lemma}
  \theoremstyle{remark}
  \newtheorem{rem}[thm]{Remark}
  \theoremstyle{plain}
  \newtheorem{prop}[thm]{Proposition}
  \theoremstyle{plain}

\newcounter{maint}

\newtheorem{mainthm}[maint]{Theorem}

 \newtheorem{step-u352}{Step}
  \newtheorem{step-l322}{Step}
\newtheorem{step-l332}{Step}
  \newtheorem{step-u33}{Step}
  \newtheorem{step-2A8}{Step}
  \newtheorem{step-A5xD5}{Step}

  \newtheorem{step-hs}{Step}
  \newtheorem{step-j2}{Step}
  \newtheorem{step-co1}{Step}
  \newtheorem{step-co2}{Step}
  \newtheorem{step-co3}{Step}
  \newtheorem{step-mcl}{Step}
  \newtheorem{step-suz}{Step}
  \newtheorem{step-j1}{Step}
  \newtheorem{step-j3}{Step}
  \newtheorem{step-ru}{Step}
  \newtheorem{step-on}{Step}
  \newtheorem{step-ly}{Step}
  \newtheorem{step-j4}{Step}
  \newtheorem{step-he}{Step}
  \newtheorem{step-hn}{Step}
  \newtheorem{step-th}{Step}
  \newtheorem{step-fi22}{Step}
  \newtheorem{step-fi23}{Step}
  \newtheorem{step-fii24}{Step}
  \newtheorem{step-b}{Step}
  \newtheorem{step-m}{Step}
  \newtheorem{step-tits}{Step}

\newcounter{algo}

\newtheorem{alg}[algo]{Algorithm}

  \theoremstyle{remark}

\newtheorem*{acknowledgement*}{Acknowledgement}
\newtheorem*{acknowledgements*}{Acknowledgements}

\theoremstyle{definition}
\newtheorem{definition}[thm]{Definition}
\newtheorem{defi}{Definition}

\theoremstyle{definition}
\newtheorem{exa}[thm]{Example}

\makeatother

\newcommand{\Z}{{\mathbb Z}}

\newcommand{\trid}{\triangleright}
\newcommand\toba{{\mathfrak B }}

\newcommand{\M}{{\mathcal M}}
\newcommand{\Q}{{\mathcal Q}}
\newcommand{\F}{{\mathbb F}}

\newcommand{\D}{{\mathcal D}}

\def\pf{\begin{proof}}

\def\epf{\end{proof}}

\begin{document}

\title[Nichols algebras over the sporadic simple groups]{Pointed Hopf algebras over the sporadic simple groups}

\author[Andruskiewitsch, Fantino, Gra\~na, Vendramin]{N.
Andruskiewitsch, F. Fantino, M. Gra\~na and L. Vendramin}

\thanks{This work was partially supported by ANPCyT-Foncyt, CONICET, Ministerio de Ciencia y
Tecnolog\'{\i}a (C\'ordoba), Secyt-UNC and Secyt-UBA}

\address{\noindent N. A., F. F. : Facultad de Matem\'atica, Astronom\'{\i}a y F\'{\i}sica,
Universidad Nacional de C\'ordoba. CIEM -- CONICET. 
Medina Allende s/n (5000) Ciudad Universitaria, C\'ordoba,
Argentina}
\address{\noindent M. G., L. V. : Departamento de Matem\'atica -- FCEyN,
Universidad de Buenos Aires, Pab. I -- Ciudad Universitaria (1428)
Buenos Aires -- Argentina}
\address{\noindent L. V. : Instituto de Ciencias, Universidad de Gral. Sarmiento, J.M. Gutierrez
1150, Los Polvorines (1653), Buenos Aires -- Argentina  }

\address{}

\email{(andrus, fantino)@famaf.unc.edu.ar} \email{(matiasg,
lvendramin)@dm.uba.ar}

\date{\today}

\begin{abstract}
We show that every finite-dimensional complex pointed Hopf algebra with group of
group-likes isomorphic to a sporadic group is a group algebra, except for the Fischer group $Fi_{22}$,  the Baby
Monster and the Monster. For these three groups, we give a short list of irreducible Yetter-Drinfeld modules
whose Nichols algebra is not known to be \fd{}.
\end{abstract}
\maketitle


\section*{Introduction}

This paper contributes to the classification of finite-dimensional
pointed Hopf algebras over an algebraically closed field $\ku$ of characteristic 0. There are different possible
approaches to this problem; one of them is to fix a finite group
$G$ and to address the classification of finite-dimensional pointed Hopf algebras $H$
such that $G(H)\simeq G$. Due to the intrinsic difficulty of this
problem, it is natural to consider separately different classes of
finite groups. A considerable progress in the case when $G$ is
abelian was achieved in \cite{AS-ann}. There are reasons to
consider next the class of finite simple, or close to simple,
groups. See \cite{AFGV} and references therein for $G$ symmetric
or alternating. Nichols algebras over Mathieu groups are partially studied
in \cite{fantino-2007}. In this paper we deal with sporadic simple
groups.

\begin{defi}\label{def-intro:gpo-collapses}
    We shall say that a finite group $G$ \emph{collapses} if
for any finite-dimensional pointed Hopf algebra $H$, with $G(H) \simeq G$, then $H\simeq \ku G$. \end{defi}

\begin{mainthm}\label{teor:complete}
If $G$ is a sporadic simple group, then it collapses, except for
the groups $G = Fi_{22}$,  $B$, $M$. For these groups,
the list of irreducible Yetter-Drinfeld modules $M(\oc, \rho)$
whose Nichols algebra is not known to be \fd{} appears in Table
\ref{tab:1}.\end{mainthm}

Notice that the Nichols algebras of \emph{reducible} Yetter-Drinfeld modules
over the groups in Table  \ref{tab:1} are infinite-dimensional \cite[8.3]{HS1}.
The conjugacy classes of sporadic groups are labeled as in the \textsf{ATLAS};
the notation for the representations in Table \ref{tab:1} is discussed in
Subsection \ref{subsect:proofs-a}.

\begin{table}[ht]
\caption{Nichols algebras  still open.}\label{tab:1}
\begin{tabular}{|c|p{1,7cm}|p{2,5cm}|c|p{1,5cm}|}
\hline $G$ & {\bf Class} & {\bf Obs.}& \small{{\bf Centralizer}} & {\bf
Rep.}
\\ \hline
\hline $Fi_{22}$  &  \textup{22A, 22B} & \small{quasi-real} $j=3$,
$g^{9}\neq g$ & $\Z/22$ &  $\chi_{-1}$
\\ \hline
\hline $B$  & \textup{16C, 16D} & real & Order 2048  &

\\ \cline{2-5}&  \textup{32A, 32B} & real & Order 128 &

\\ \cline{2-5}&  \textup{32C, 32D} & \small{quasi-real} $j=3$,
$g^{9}\neq g$ & Order 128 &

\\ \cline{2-5}&  \textup{34A} & real &  $\Z/34\times \Z/2$ & $\chi_{-1} \otimes \epsilon$, \newline $\chi_{-1} \otimes \sgn$

\\ \cline{2-5}&  \textup{46A, 46B} & \small{quasi-real} $j=3$,
    $g^{9}\neq g$ & $\Z/46$ & $\chi_{-1}$

\\ \hline \hline $M$  &
\textup{32A, 32B} & real  & Order 128 &

\\ \cline{2-5}& \textup{46A, 46B} & \small{quasi-real} $j=3$,
$g^{9}\neq g$  & $\Z/23 \times \mathbb D_4$ & $\epsilon \otimes
\rho_2$


\\ \cline{2-5}& \textup{92A, 92B} & \small{quasi-real} $j=3$,
$g^{9}\neq g$  & $\Z/92$ & $\chi_{-1}$

\\ \cline{2-5}&  \textup{94A, 94B} & \small{quasi-real} $j=3$,
$g^{9}\neq g$ & $\Z/94$& $\chi_{-1}$
  \\ \hline
\end{tabular}
\end{table}

\bigbreak Another approach to the mentioned classification problem is
through Nichols algebras associated to racks. Precisely, one has
to attack the following question:

\begin{center}
    \emph{For every finite indecomposable rack $X$, for every $n\in \N$, and for every
    $2$-cocycle\\
    $q\in Z^2(X, \GL(n, \ku^\times))$, determine if $\dim \toba(X,q) < \infty$}.
\end{center}

Here $\toba(X,q)$ denotes the Nichols algebra associated to the
pair $(X, q)$, see Subsection \ref{subsect:nichols}. We refer to
\cite{AG1}, \cite[Introduction]{AFGV} for the relation between
this question and the classification problem. Again, it is natural
to consider separately different classes of finite racks; again,
it is natural to start by the class of finite simple racks. These
were classified in \cite[Th. 3.9, Th. 3.12]{AG1}, see also [J];
among them, there are the racks arising as non-trivial conjugacy
classes of finite simple groups. A substantial part of the proof
of Theorem \ref{teor:complete} follows from a much more general
result for a large family of conjugacy classes of sporadic
groups. It is convenient to introduce the following terminology
before stating our next main theorem.

\begin{defi}\label{def-intro:collapses}
    We shall say that a finite simple rack $X$ \emph{collapses} if $\dim \toba(X,q) = \infty$
    for any $2$-cocycle $q\in Z^2(X, \GL(n, \ku^\times))$, for any $n\in \N$.
\end{defi}

See Definition \ref{def:tiposbyd} for the notion of rack of type D; such rack collapses by
\cite[Th. 3.6]{AFGV}, recalled in Theorem \ref{th:racks-claseD}. This result is a translation to
the language of racks of \cite[Th. 8.6]{HS1}, whose proof relies on  results from
\cite{AHS}.

\begin{mainthm}\label{th:racks-liquidados}
  If $G$ is a sporadic simple group and  $\oc$ is a non-trivial conjugacy class of
  $G$ NOT listed in Table \ref{tab:0}, then $\oc$ is of type D, hence it collapses.
\end{mainthm}

\begin{table}[ht]
\begin{center}
\caption{Conjugacy classes not known of type D; \newline those which are NOT of type D appear in bold.}\label{tab:0}
\begin{tabular}{|p{1cm}|c||p{1cm}|c|}
\hline $G$ & {\bf Classes} & $G$ & {\bf Classes}
\\ \hline  $M_{11}$ &  \textup{\bf 8A, 8B, 11A, 11B} &
$Co_{1}$ &  \textup{3A, 23A, 23B}
\\ \hline  $M_{12}$ &  \textup{\bf 11A, 11B} &
$J_{1}$ &  \textup{\bf 15A, 15B, 19A, 19B, 19C}
\\ \hline  $M_{22}$ &  \textup{\bf 11A, 11B} &  $O'N$ &  \textup{31A, 31B}
\\ \hline  $M_{23}$ &  \textup{\bf 23A, 23B} &
$J_{3}$ &  \textup{5A, 5B, 19A, 19B}
\\ \hline  $M_{24}$ &  \textup{\bf 23A, 23B} &    $Ru$ &  \textup{29A, 29B}
\\ \hline  $J_{2}$ &  \textup{\bf 2A, 3A} &    $He$ &  \textup{all collapse}
\\ \hline  $Suz$ &  \textup{3A} &  $Fi_{22}$ &  \textup{{\bf 2A}, 22A,
22B}
\\ \hline  $HS$ &  \textup{11A, 11B} &   $Fi_{23}$ &  \textup{{\bf 2A},  23A, 23B}
\\ \hline  $McL$ &  \textup{11A, 11B} &  $HN$ &  \textup{all collapse}
\\ \hline  $Co_{3}$ &  \textup{23A, 23B} & $Th$ & \textup{all collapse}
\\ \hline  $Co_{2}$ &  \textup{{\bf 2A}, 23A, 23B} &  $T$ &  \textup{2A}
\\ \hline
\end{tabular}
\end{center}
\begin{tabular}{|p{0,8cm}|c|}
\hline $G$ & {\bf Classes}
\\ \hline  $Ly$ &  \textup{33A, 33B, 37A,
37B, 67A, 67B, 67C}
\\ \hline   $J_4$ &  \textup{29A, 37A, 37B, 37C, 43A, 43B, 43C}
\\ \hline   $Fi'_{24}$ &  \textup{23A, 23B, 27B, 27C, 29A, 29B, 33A, 33B, 39C, 39D}
\\ \hline $B$ &  \textup{2A, 16C, 16D, 32A, 32B, 32C, 32D,
34A, } \\ \hline & \textup{46A, 46B, 47A, 47B}
\\ \hline $M$ &  \textup{32A, 32B, 41A,
46A, 46B, 47A, 47B, 59A, 59B, } \\ \hline & \textup{69A, 69B, 71A, 71B, 87A, 87B, 92A, 92B, 94A, 94B} 
\\ \hline\end{tabular}
\end{table}

\bigbreak Theorem \ref{th:racks-liquidados} is not merely an auxiliary step
towards Theorem \ref{teor:complete}; its consequences, illustrated
by Theorems \ref{th:racks-liquidados-aplicacion} and
\ref{th:factores-colapsan}, show how crucial is the use of racks
in the classification of pointed Hopf algebras.

\bigbreak Our proofs of Theorems \ref{teor:complete}, \ref{th:racks-liquidados}
and \ref{th:racks-liquidados-aplicacion} are based on reductions to problems on
conjugacy classes of finite groups; we solve these problems in the present
setting with the help of \cite{GAP}. For completeness, we review these
reductions in Section \ref{sect:preliminaries}, and discuss the elements of
\textsf{GAP} that we are using. We discuss in Section \ref{sect:proofs-b} the
algorithms used to justify Theorem \ref{th:racks-liquidados}.  We discuss
Theorem \ref{teor:complete} in Section \ref{sect:proofs-a}. We treat the
classes not covered by Theorems \ref{th:racks-liquidados} in Subsection
\ref{subsect:proofs-a}.  In Subsection \ref{subsect:remarks-remaining}, we
explain why the remaining classes in Table \ref{tab:1} are beyond our present
knowledge.

\bigbreak In this paper we explain the algorithms, and their theoretical support, for the proofs of our main results.
The details of how these algorithms are actually implemented are given in the companion paper
\cite{logbook}. We believe, and hope, that these details are enough to guide the diligent reader to repeat and corroborate
our calculations.

\section{Preliminaries}\label{sect:preliminaries}

\subsection{Nichols algebras}\label{subsect:nichols} A \emph{braided vector space} is a pair $(V, c)$,
where $V$ is a vector space and $c\in\mathbf{GL}(V\otimes V)$ is a
solution of the braid equation, that is
\begin{align}\label{eqn:braid-eq}
(c\otimes \id)(\id\otimes c)(c\otimes \id) =  (\id\otimes
c)(c\otimes \id)(\id\otimes c).
\end{align}
The  braid equation is equivalent to the celebrated quantum
Yang-Baxter equation, that plays an important role in statistical
mechanics.

There is a very interesting object associated to a braided vector
space, \emph{its Nichols algebra}, defined as follows.
\begin{itemize}
    \item The solution $c$ of the braid equation \eqref{eqn:braid-eq} induces a representation of the braid group
    $\bn$ in the  $n$-th tensor product $V^{\otimes n}$.

    \item Let $M: \sn \to \bn$ be the so-called Matsumoto section; it preserves the product
    only when the length is preserved.

    \item Let $\Omega_n = \sum_{\sigma \in \sn} M(\sigma) \in \End (V^{\otimes n})$ and let
    $J = \oplus_{n\ge 2} \ker \Omega_n$.    The Nichols algebra of the braided vector space $(V,c)$ is
    $\toba(V,c) = T(V)/J$.
        \end{itemize}

The study of Nichols algebras is unavoidable in the classification problem
of pointed Hopf algebras. A notable example of a Nichols algebra is the positive
part $U^+_q(\mathfrak{g})$ of a quantized enveloping algebra when $q$ is not a root of 1 \cite{Lu, Ro, Sbg}.
It is in general very difficult to compute explicitly the Nichols algebra of
a braided vector space -- that is, to determine its dimension, or even better, an efficient set of generators.
We refer to \cite{AS-cambr} for a detailed discussion of alternative definitions,
basic results and examples of Nichols algebras. In the present paper, we shall
make use of the following observation:

\centerline{\emph{If $(W, c)$ is a braided vector subspace of $(V, c)$, then
$\toba(W, c) \hookrightarrow\toba(V, c)$.}}

In particular, if we find a braided vector subspace $(W, c)$ of
$(V, c)$ with infinite-dimensional Nichols algebra, then $\dim
\toba(V, c) = \infty$ too. We shall reduce the search of a
suitable braided vector subspace to problems in finite group
theory, and then we will solve them for sporadic groups using
\textsf{GAP}.

\subsection{Yetter-Drinfeld modules}\label{subsect:yd} We now describe the class of braided
vector spaces whose Nichols algebras we need to show that have
infinite dimension. Let $G$ be a finite group. We denote by $Z(G)$
the center of $G$ and by $\Irr G$ the set of isomorphism classes
of irreducible representations of $G$. If $g\in G$, we denote
by $C_{G}(g)$ the centralizer of $g$ in $G$. The conjugacy class
of $g$ is denoted by $\oc_g$ or by $\oc^G_g$, if emphasis on the group
is needed.

\smallbreak A \emph{Yetter-Drinfeld module} over $G$ is a $G$-graded vector
space $M = \oplus_{g\in G} M_g$ provided with a $G$-module structure such that
$g\cdot M_t = M_{gtg^{-1}}$ for any $g,t\in G$. The category $\ydg$ of
Yetter-Drinfeld modules over $G$ is semisimple and its irreducible objects are
all of the following form. Let $\oc$ be a conjugacy class of $G$, $g\in \oc$
fixed, $\rho\in \Irr C_{G}(g)$. We describe the corresponding irreducible
Yetter-Drinfeld module $M(\oc, \rho)$. Let $g_1 = g$, \dots, $g_{m}$ be a
numeration of $\oc$ and let $x_i\in G$ such that $x_i \trid g = g_i$ for all
$1\le i \le m$. Then
\[
M(\oc, \rho) =\operatorname{Ind}^G_{C_G(g)}V=\oplus_{1\le i \le m} x_i\otimes V.
\]
Let $x_iv := x_i\otimes v \in M(\oc,\rho)$, $1\le i \le m$, $v\in V$.
The Yetter-Drinfeld module $M(\oc,\rho)$
is a braided vector space with braiding given by
\begin{equation}
\label{yd-braiding}
c(x_iv\otimes x_jw) = g_i\cdot(x_jw)\otimes x_iv =
x_h\,\rho(\gamma)(w) \otimes x_iv
\end{equation}
for any $1\le i,j\le m$, $v,w\in V$, where $g_ix_j = x_h\gamma$ for unique $h$,
$1\le h \le m$ and $\gamma \in C_{G}(g)$. The Nichols algebra\footnote{We omit to mention the braiding $c$ in the notation of a
Nichols algebra from now on.} of $M(\oc,\rho)$ is simply denoted $\toba(\oc,\rho)$.

\smallbreak Let $H$ be a pointed Hopf algebra with $G(H)\simeq G$. Then there are
two fundamental invariants of $H$, a Yetter-Drinfeld module $V$ over $\ku G$
(called the infinitesimal braiding of $H$) and its Nichols algebra
$\toba(V)$, see \cite{AS-cambr,A}.  A basic problem for the
classification of \fd{} pointed Hopf algebras over $G$ is the determination of
all  Yetter-Drinfeld modules $V$ over $\ku G$ such that the Nichols algebra
$\toba(V)$ is \fd{}. In particular, the following statements are equivalent:

\smallbreak
 \begin{asparaenum}[(1)]
     \item \emph{If $H$ is a \fd{} pointed Hopf
         algebra with $G(H) \simeq G$, then $H\simeq \ku G$.}

        \smallbreak
     \item \emph{If $V\neq 0$ is a Yetter-Drinfeld module over $\ku G$,
         then $\dim \toba(V) = \infty$.}

        \smallbreak \item \emph{If $V$ is an \emph{irreducible} Yetter-Drinfeld module over $\ku G$,
         then $\dim \toba(V) = \infty$.}\label{it:viejo3}
 \end{asparaenum}

\smallbreak Therefore, for a fixed group $G$, we aim to know when
$\dim \toba(V) = \infty$ for an irreducible $V\in \ydg$; that is, when $\dim \toba(\Oc, \rho) = \infty$ for a pair $(\Oc, \rho)$ as above. As we
said, we shall look at braided vector subspaces of $V$, and to
describe them we find convenient the language of racks and
cocycles.

\subsection{Racks}\label{subsect:racks} A \emph{rack} is a pair $(X,\trid)$ where $X$ is a
non-empty set and $\trid:X\times X\to X$ is a function such that
$\phi_i:X\to X$, $\phi_i (j) := i\trid j$, is a bijection for all
$i\in X$, and $i\trid(j\trid k)=(i\trid j)\trid(i\trid k)$ for all
$i$, $j$, $k\in X$. For instance, a group $G$, and any conjugacy
class in $G$, is a rack with $x\trid y=x y x^{-1}$. In this case,
$j\trid i=i$ whenever $i\trid j=j$ and $i\trid i = i$ for all
$i\in G$. A rack $X$ is \emph{abelian} if for any $x,y\in X$,
$x\trid y = y$. See \cite{AG1} for a survey on racks.

\smallbreak
Let $X$ be a rack. Given $n\in \N$, a map $q:X\times
X\to\GL(n,\ku)$ is a \emph{2-cocycle} if $$q_{x,y\trid z}q_{y,z}=
q_{x\trid y,x\trid z}q_{x,z},$$ for all $x,y,z\in X$. Let $q$ be a
2-cocycle, $V= \ku X\otimes\ku^{n}$, where $\ku X$ is the vector
space with basis $e_x$, for $x\in X$. We denote $e_xv := e_x\otimes v$.
Consider the linear isomorphism
$c^q:V\otimes V\to V\otimes V$,
$$c^q(e_xv\otimes e_yw)=e_{x\trid y}q_{x,y}(w)\otimes e_xv,$$
$x\in X$, $y\in X$, $v\in\ku^{n}$, $w\in\ku^{n}$. Then $c^q$ is a
solution of the braid equation; its Nichols algebra is denoted
$\toba(X,q)$. Pointed Hopf algebras are related to Nichols
algebras over racks by the following observations, see \cite[Th.
4.14]{AG1}.

\begin{itemize}
\item If $X$ is a rack, $n\in \N$,
and $q:X\times X\to\GL(n,\ku)$ is a 2-cocycle, then there exists a
group $G$ such that $V = \ku X\otimes \ku^n$ is a Yetter--Drinfeld
module over $G$ and the braiding of $V$ as an object in $\ydg$
coincides with $c^q$. If the image of $q$ is contained in a finite
subgroup $\varGamma\subset\GL(n,\ku)$, then $G$ can be chosen to
be finite.

\medbreak\item Conversely, if $G$ is a finite group and $V =
M(\oc, \rho)\in\ydg$ is irreducible, then there exists a
finite subgroup $\varGamma$ of $\GL(n,\ku)$, $n = \dim \rho$, and a 2-cocycle
$q:X\times X\to\varGamma$ such that $V$ is given as above and the
braiding $c\in\Aut(V\otimes V)$ in the category $\ydg$ coincides
with $c^q$.
\end{itemize}

\subsection{Abelian techniques}
Let $G$ be a finite group, $\oc$ a conjugacy class of $G$, $g\in
\oc$  and $(\rho,V)\in \Irr C_{G}(g)$. Our next goal is to
describe techniques to conclude that $\dim\mathfrak{B}(\oc,\rho) =
\infty$ from the analysis of abelian subracks; these will be
applied in the proof of Theorem \ref{teor:complete}. Abelian subracks give rise
to braided subspaces of $M(\Oc,\rho)$ of diagonal type; then the classification
of braided  vector spaces of diagonal type with finite-dimensional Nichols algebra \cite{H-all} may be invoked.
We shall freely use the notations and results from \cite{H-all}, in particular the notion of \emph{generalized Dynkin diagram}.
We start by subracks with one element. By the Schur lemma, $\rho(g)$ is a
scalar. It is well-known that
\begin{equation}\label{eqn:rho-g-uno}
\rho(g) = 1 \implies \dim\mathfrak{B}(\oc,\rho)=\infty.
\end{equation}

\medbreak We consider next subracks with two or three elements.

\begin{lem}
\label{lem:inversos} \cite[2.2]{AZ} Assume that
$\dim\mathfrak{B}(\oc,\rho)<\infty$. If $g\in G$ is
\emph{real}\index{real}, that is $g^{-1}\in\oc$, then
$\rho(g)=-1$. In particular, the order of $g$ is even. \qed
\end{lem}

If $g$ is not real, but it is conjugated to $g^{j}\neq g$ for some
$j\in \Z$, then we shall say that $g\in G$ and $\oc$ are
\emph{quasi-real}.\index{quasi-real}
\begin{lem} \cite[1.8 and 1.9]{AF2}, \cite[2.2]{FGV}
\label{lem:potencias}Assume that
$\dim\mathfrak{B}(\oc,\rho)<\infty$ and that there exists $j$ such
that $g\ne g^{j}\in\oc$.

\renewcommand{\theenumi}{\roman{enumi}}\renewcommand{\labelenumi}{(\theenumi)}
\begin{enumerate}
\item If $\deg\rho>1$, then $\rho(g)=-1$ and $g$ has even order.

\item If $\deg\rho=1$, then $\rho(g)=-1$ and $g$ has even order or $\rho(g)\in\mathcal{R}_{3}$.

\item If $g^{j^{2}}\ne g$, then $\rho(g)=-1$. \qed
\end{enumerate}
\end{lem}


\bigbreak Involutions, that is elements of order 2, are real but Lemma \ref{lem:inversos} does not give useful information.
Our next criterium deals with classes of involutions and representations of degree greater
than one. Another useful criterium for classes of involutions is Proposition \ref{cor:ct_a4} below.

\begin{lem}\label{lem:contar}
Let $G$ be a finite group, $\oc$ the conjugacy class of an
involution $g\in G$ and $(\rho,V)\in\Irr C_G(g)$. Assume that
there exists an involution $x$ such that $h=xgx$ and $gh=hg$. If
$\rho(h)$ has $a$ eigenvalues $1$, $b$ eigenvalues $-1$, where
either $a>0$ and $b>3$ or $a>3$ and $b>0$, then
$\dim\mathfrak{B}(\mathcal{O},\rho)=\infty$.
\end{lem}

\begin{proof}
By \eqref{eqn:rho-g-uno} we can assume $\rho(g)=-1$. Let $x_1= e$
and $x_2=x$. Since $gh=hg$, there exists a linear basis
$\{v_1,\ldots,v_n\}$ of $V$ such that $\rho(g)$ and $\rho(h)$ are
simultaneously diagonalizable in this basis. Define $W$ as the
subspace generated by $x_{\ell} v_i=x_{\ell}\otimes v_i$, where
$\ell=1$, $2$ and $1\leq i\leq n$.  Then $W$ is a braided vector
space of diagonal type with braiding given by $c(x_{\ell}
v_i\otimes x_{\ell} v_j)=-x_{\ell} v_j\otimes x_{\ell} v_i$ and
\begin{gather*}
c(x_1 v_i \otimes x_2 v_j)= x_2\rho(h) v_j\otimes x_1 v_i,\\
c(x_2 v_i \otimes x_1 v_j)= x_1\rho(h) v_j\otimes x_2 v_i,
\end{gather*}
for $\ell=1$, $2$, and  $1\leq i, j\leq n$. Therefore, the results
follows from \cite{H-all} because the generalized Dynkin diagram has at least
one vertex with valency $>3$. For instance assume that $a=4$ and $b=1$. 
Then the generalized Dynkin diagram is 

\begin{figure}[ht]
\caption{\label{fi:triangulito_a=4b=1}}
\begin{center}
\vspace{1cm}
\begin{align*}
\setlength{\unitlength}{1.4cm}
\begin{picture}(1,0)
\put(-1,0){\circle*{.15}} \put(0,0){\circle*{.15}}\put(1,0){\circle*{.15}}\put(2,0){\circle*{.15}}\put(3,0){\circle*{.15}}
\put(-1,1){\circle*{.15}}\put(0,1){\circle*{.15}} \put(1,1){\circle*{.15}}\put(2,1){\circle*{.15}}\put(3,1){\circle*{.15}}
\put(-1,0){\line(4,1){4}}
\put(0,0){\line(3,1){3}}\put(1,0){\line(2,1){2}}\put(2,0){\line(1,1){1}}
\put(-1,1){\line(4,-1){4}}  \put(0,1){\line(3,-1){3}} \put(1,1){\line(2,-1){2}} \put(2,1){\line(1,-1){1}}
\end{picture}\qquad \qquad
\end{align*}
\end{center}
\end{figure}
This completes the proof.                                                  
\end{proof}

The computation of the multiplicities $a$ and $b$ in the statement
of the lemma can be performed using the following remark.

    \begin{rem}\label{rem:contar_autovalores}
        Let $G$ be a finite group, $g\in G$ and $(\rho,V)$ a
        representation of $G$. The multiplicities of the eigenvalues of
        $\rho(g)$ are the scalar products of the restriction of the
        representation to the cyclic group generated by $g$ with the
        irreducible characters of this cyclic group. The \textsf{GAP} function
        \textsf{EigenvaluesChar} can be used for this computation.
    \end{rem}
%

\subsection{The technique of the subgroup}
Let $G$ be a finite group, $\sigma\in G$, $\oc^G_\sigma = \oc^G$
its conjugacy class, $C_G(\sigma)$ its centralizer and $\rho\in
\Irr C_G(\sigma)$. If $H$ is a subgroup of $G$ and $\sigma\in H$,
then $\Oc_{\sigma}^{H} = \oc^H$ denotes the conjugacy class of
$\sigma$ in $H$. Let $\rho\vert_{C_H(\sigma)} =
\tau_1\oplus\dots\oplus\tau_s$ where $\tau_j\in \Irr C_H(\sigma)$,
$1\le j \le s$.

\begin{lem}\label{lema:technique-fusionsubgroups}\label{lem:subgrupo_general}  \cite[3.2]{AFGV}
    \renewcommand{\theenumi}{\roman{enumi}}   \renewcommand{\labelenumi}{(\theenumi)}
    \begin{enumerate}
        \item\label{item:subgpo1}
            If $\dim \toba(\oc^H,\lambda)=\infty$ for all
            $\lambda\in\Irr C_{H}(\sigma)$, then
            $\dim\toba(\oc^G,\rho)=\infty$ for all $\rho\in
            \Irr C_{G}(\sigma)$.

        \medbreak
        \item\label{item:subgpo3}
            Let $\sigma_1,\sigma_2\in\oc^G\cap H$. Let
            $\mathcal{O}_i=\mathcal{O}_{\sigma_i}^H$ and
            assume that $\mathcal{O}_1\ne\mathcal{O}_2$. If
            $\dim\mathfrak{B}(M(\mathcal{O}_{1},\lambda_{1})\oplus
            M(\mathcal{O}_{2},\lambda_{2}))=\infty$ for all
            pairs $\lambda_1\in\Irr C_H(\sigma_1)$,
            $\lambda_2\in\Irr C_H(\sigma_2)$, then
            $\dim\mathfrak{B}(\mathcal{O}^{G},\rho)=\infty$
            for all $\rho\in \Irr C_{G}(\sigma)$. \qed

    \end{enumerate}
\end{lem}

\subsection{The group $\ac$ and conjugacy classes of involutions}

In this section we give another criterium for Nichols algebras
over the conjugacy class of an involution.

\begin{lem}
Let $g,h\in G$ with $g$ an involution, $h$ and $gh$ of order $3$;
let $\oc$ be the conjugacy class of $g$. Then
$\dim\mathfrak{B}(\oc, \rho)=\infty$.
\end{lem}\label{lem:a4}
\begin{proof}
The alternating group in four letters $\mathbb{A}_{4}$ can be
presented by generators $g$ and $h$ with relations $g^{2}=
h^{3}=(gh)^{3}= e$ \cite{dickson}. Thus, the subgroup $H$ of $G$
generated by $g$ and $h$ is isomorphic to $\ac$. Indeed, the
hypothesis implies that the elements $g,h,h^2, gh, (gh)^2$ of $G$
are all distinct, but a proper subgroup of $\ac$ has at most 4
elements. From \cite{AF2}, the conjugacy class of involutions in
$\mathbb{A}_{4}$ gives infinite-dimensional Nichols algebras for
every representation. Hence Lemma
\ref{lema:technique-fusionsubgroups} applies.
\end{proof}

Now we need an efficient way of checking if $\mathbb{A}_{4}$ is a
subgroup. To this purpose, first recall a very useful result from
the theory of groups.

\begin{prop}
\label{prop:go}
\cite[Th. 4.2.12]{g}. Let $G$ be a finite group and let $\oc_{i}$, $\oc_{j}$, $\oc_{k}$
be some conjugacy classes. If
$S(\oc_{i},\oc_{j},\oc_{k})$ is the number of times that a given
element of $\oc_{k}$ can be expressed as an ordered product of an
element of $\oc_{i}$ with an element of $\oc_{j}$, then
\begin{equation}\label{eqn:suma-a4}
S(\oc_{i},\oc_{j},\oc_{k})=\frac{|\oc_{i}||\oc_{j}|}{|G|}
\sum_{\chi}\frac{\chi(\oc_{i})\chi(\oc_{j})\overline{\chi(\oc_{k})}}{\chi(1)},
\end{equation}
where $\chi$ runs over all irreducible characters of $G$. \qed
\end{prop}

The last proposition will be used in connection with the following criterium.

\begin{prop}
\label{cor:ct_a4} Let $\oc$ be the conjugacy class of an
involution, and $\mathcal{K}_{1},\mathcal{K}_{2}$ conjugacy
classes of elements of order $3$. If
$S(\oc,\mathcal{K}_{1},\mathcal{K}_{2})\geq 1$ then
$\dim\mathfrak{B}(\oc, \rho)=\infty$.
\end{prop}
\begin{proof}
By hypothesis, there exists $a\in\oc$ and
$b_{i}\in\mathcal{K}_{i}$, $i=1,2$, such that $ab_{1}=b_{2}$. Then
Lemma \ref{lem:a4} applies.
\end{proof}

\medbreak\subsection{A technique with simple affine racks}\label{subsec:hs-conj}
We now apply \cite[Th. 8.2]{HS1} to racks arising as union of two affine racks.
This is used for the conjugacy classes labeled  \textup{8A, 8B} of the Mathieu group $M_{11}$,
that can not be treated otherwise.

If $G$ is a finite group and $g,h \in G$, then the
conjugacy classes $\Oc_g$ and $\Oc_h$ \emph{commute} if $st = ts$
for all $s\in \Oc_g$, and for all $t\in \Oc_h$. Let
$$
\Fc(G) = \{\Oc \text{ conjugacy class of }G: \dim \toba (\Oc,
\rho) < \infty \text{ for some } \rho\}.
$$

\begin{thm}\label{thm:hs2} \cite[Th. 8.2]{HS1}
    Let $G$ be a finite group such that any two conjugacy classes in $\Fc(G)$ do not commute.
    Let $0\neq U\in \ydg$  such that $\dim \toba(U) < \infty$. Then $U$ is irreducible. \qed
\end{thm}

If $A$ is an abelian group and $T\in\Aut(A)$, then $A$ becomes a rack with $x\trid
y=(1-T)x+Ty$. It will be denoted by $(A,T)$ and called an \emph{affine
rack}. We realize it as a conjugacy class in
the semidirect product $G=A\rtimes\langle T\rangle$. The conjugation in $G$ gives
\begin{equation}\label{eqn:affinerackmult}
    (v, T^{h})\trid (w, T^j) = (T^{h}(w) + (\id - T^j)(v) , T^j).
\end{equation}
We denote $\Q_{A, T}^{j} := \{(w, T^j): w\in A\}$, $j\in \Z/d$, a
subrack of $G$ isomorphic to the affine rack $(A, T^{j})$.

\bigbreak We assume that $(A,T)$ is a {\em simple} affine rack; that is, $A =
\F_{\hspace{-2pt}p}^{\, t}$, $p$ a prime, and $T\in \GL(t, \F_{\hspace{-2pt}p})
- \{\id\}$ of order $d$, acting irreducibly.

\begin{lem}\label{lema:dobleafin-82} Suppose that $p >2$ when $d$ is even. Let $G = A \rtimes \langle
T\rangle$. If $0\neq U\in \ydg$ satisfies $\dim \toba(U) < \infty$, then
$U$ is irreducible.
\end{lem}

\pf

By Theorem \ref{thm:hs2}, we have to show that any two conjugacy
classes in $\Fc(G)$ do not commute. We claim that
\renewcommand{\theenumi}{\alph{enumi}}\renewcommand{\labelenumi}{(\theenumi)}
\begin{enumerate}
        \item\label{clases-semidirecto} the conjugacy
classes of $G$ are either $\Q_{A, T}^{j}$ with $j\neq 0$, or else
the orbits of $T$ in $A$.

        \item\label{clases-semidirecto-nichols-infinita}
        $\Fc(G)\subset \{\Q_{A, T}^{j}: \, j\neq 0\}$; hence any two conjugacy
classes in $\Fc(G)$ do not commute, cf.
\eqref{eqn:affinerackmult}.
\end{enumerate}

Part (a) is elementary, but we sketch the argument. It is evident from
\eqref{eqn:affinerackmult} that the conjugacy class of $(w,\id)$ is the orbit
of $w$ under $T$. If $0 < j < d$, then $\id - T^j$ is bijective, since its
kernel is a $T$-invariant subspace, but we are assuming that $T$ is
irreducible. Hence the class of $(0, T^j)$ is $\Q_{A, T}^{j}$.

We prove (b). Let $v\in A -0$; the centralizer of $v$ is $A$. Set
$\sigma_k = T^k(v)$ and $x_k = (0, T^k)$; thus $x_k\trid \sigma_0
= \sigma_k$ and $\sigma_{\ell}x_k = x_k\sigma_{\ell - k}$, $0\le
\ell, k \le d-1$. Let $\chi \in \Irr A$; then the braiding in
$M(\mathcal O_v, \chi)$ is given by $c(x_k\otimes x_{\ell}) =
\chi(\sigma_{k - \ell}) x_{\ell} \otimes x_k$. In other words,
this is of diagonal type with matrix $q_{k\ell} = \chi(T^{k -
\ell}(v))$. Let $\Delta$ be the generalized Dynkin diagram
associated to $(q_{k\ell})$. Now we can identify $A = \F_q$, with
$q = p^{t}$, in such a way that $T(v) = \xi v$, where $\xi \in
\F_q^{\times}$ has order $d$. Hence $q_{k\ell} q_{\ell k} =
\chi\big((\xi^{\ell - k} + \xi^{k - \ell})v\big)$.   Notice that
the order of $q_{k\ell} q_{\ell k}$ divides $p$. Also,
\begin{equation}\label{eqn:xi}
1+\xi+\xi^2+\cdots+\xi^{d-1}=0.
\end{equation}

We can assume $\chi(v)\neq 1$. Now, we consider different cases.

Suppose that $d$ is odd.
\begin{itemize}
\item[(I)] Assume that $3 \nmid d$. If there exists $\ell$, with $1\leq \ell \leq d-1$,
such that $q_{0\ell}q_{\ell\, 0}\neq 1$, then $\dim\toba(\mathcal O_v, \chi)=\infty$.
Indeed, $\Delta$ contains a cycle of length greater than 3, namely the set of vertices $0$, $\ell$, $2\ell$,
\dots, $(M-1)\ell$, where $M$ is the order of $\ell$. Then the result follows from \cite{H-all}.
If there is not such $\ell$, then $ v = - \sum_{1\le \ell\le (d-1)/2} (\xi^{\ell} + \xi^{d - \ell})v \in \ker \chi$
by \eqref{eqn:xi}, a contradiction.

\medbreak\item[(II)] Suppose that $3 \mid d$.
\begin{itemize}
  \item[(i)] If there exists $\ell$, with $1\leq \ell \leq d-1$ and $\ell \neq \frac{d}{3}, \frac{2d}{3}$, such that
  $q_{0\ell}q_{\ell \, 0}\neq 1$, then $\dim\toba(\mathcal O_v, \chi)=\infty$ because $\Delta$ contains a cycle as in (I).

\medbreak\item[(ii)] Assume that $q_{0\ell}q_{\ell \, 0}= 1$ for all $\ell$, with $1\leq \ell \leq d-1$ and $\ell \neq \frac{d}{3}, \frac{2d}{3}$.
  If $q_{0\frac{d}{3}}q_{\frac{d}{3} 0} = q_{0\frac{2d}{3}}q_{\frac{2d}{3} 0} = 1$,
  then $\chi(v)=1$ by \eqref{eqn:xi}, a contradiction. On the other hand, if $\omega := q_{0\frac{d}{3}}q_{\frac{d}{3} 0}
  =q_{0\frac{2d}{3}}q_{\frac{2d}{3} 0}\neq 1$, then $\Delta$ contains a triangle like in Figure \ref{fi:triangulito},
  but no triangle of this sort appears in \cite[Table 2]{H-all}.
\end{itemize}
\end{itemize}

\begin{figure}[ht]
\caption{}\label{fi:triangulito}
    \vspace{1cm}
    \begin{align*}
        \setlength{\unitlength}{1.4cm}
        \begin{picture}(1,0)
            \put(0,0){\circle*{.15}} \put(1,1){\circle*{.15}}
            \put(2,0){\circle*{.15}} \put(0,0){\line(1,1){1}}
            \put(0,0){\line(2,0){2}} \put(1,1){\line(1,-1){1}}
            \put(0,.7){$\omega$}
            \put(1.6,.7){$\omega$} \put(0.9,-0.4){$\omega$}
            \put(-0.6,-0.07){$\chi(v)$} \put(2.2,-0.07){$\chi(v)$}
            \put(0.9,1.2){$\chi(v)$}
        \end{picture}\qquad \qquad
    \end{align*}
\end{figure}

Assume that $d$ is even. Then $q_{0\frac{d}{2}}q_{\frac{d}{2} 0}= \chi(v)^{-2} \neq 1$ because $p>2$; hence,
$0$ is connected to $\frac{d}{2}$ and the sub-diagram spanned by 0
and $\frac{d}{2}$ is of Cartan type $A_1^{(1)}$. By \cite{H-all}, $\dim\toba(\mathcal O_v, \chi) = \infty$.
In conclusion, $\mathcal O_v\notin \Fc(G)$. This proves (b). \epf

\begin{prop}\label{lema:dobleafin-8.2}
Let $G$ be a finite group;  $(A,T)$ an affine simple rack with
$\vert T\vert = d$ and $p > 2$ as above; and $\psi:A\rtimes
\langle T\rangle\to G$ a monomorphism of groups. Assume that the
conjugacy class $\oc$ of $\sigma = \psi(0,T)$ is quasi-real of
type $j$, $1 < j < d$. If $\rho\in \Irr C_{G}(\sigma)$, then $\dim
\toba(\Oc, \rho) = \infty$.
\end{prop}

\pf
This follows from Lemmata \ref{lem:subgrupo_general}
\eqref{item:subgpo3} and \ref{lema:dobleafin-82} applied to $H =$ the image of $\psi$. Indeed, $\oc_{\sigma}^H$
and $\oc_{\sigma^j}^H$ are both contained in $\oc$, but $\oc_{\sigma}^H\cap\oc_{\sigma^j}^H = \emptyset$.
\epf

\begin{exa}\label{exa:affin-notipoD}
Let $T\in  \GL(2, \F_{\hspace{-2pt}3})$ of order 8 and  $G =
\F_{\hspace{-2pt}3}^{\, 2} \rtimes \langle T\rangle \simeq (\Z/3 \oplus \Z/3)
\rtimes \Z/8$. Then there is a monomorphism of groups $\psi$ from $G$ to the
Mathieu group $M_{11}$ such that $\psi(T)$ belongs to the class $\oc$ labeled
\textup{8A} (resp. \textup{8B}) in \atlas, which is quasi-real of type 3. Then
$\dim \toba(\Oc, \rho) = \infty$ for any $\rho$.
\end{exa}

\subsection{The dihedral group and conjugacy class of involutions}

Let $n$ be an even number. Recall that the dihedral group of $2n$ elements is given by
\[
\mathbb{D}_{2n}=\langle r,s\mid r^n=s^2=1,\;srs=r^{-1}\rangle.
\]
The involutions of $\mathbb{D}_{2n}$ split in three
conjugacy classes $S=\{r^{2i}s\mid 0\leq i\leq\frac{n}{2}-1\}$, 
$R=\{r^{2i+1}s\mid 0\leq i\leq\frac{n}{2}-1\}$ and $\{r^{n/2}\}$. 

\begin{lem}
\label{lem:D4n_tipoD}
Let $n>4$ be an even number.  Then $Y=R\coprod S$ is of type D.
\end{lem}

\begin{proof}
Let $\sigma_1=s\in S$ and $\sigma_2=rs\in R$; clearly,  
$(\sigma_1\sigma_2)^2\ne(\sigma_2\sigma_1)^2$.
\end{proof}

\begin{lem}
\label{lem:D_2n}
Let $\mathcal{A}$ be a conjugacy class of involutions in a finite group $G$ and
let $\mathcal{B}$ be a conjugacy class with representative of order $n$, with
$n>2$ even. If $S(\mathcal{A},\mathcal{A},\mathcal{B})>0$, then the conjugacy
class $\mathcal{A}$ is of type D.
\end{lem}

\begin{proof}
    Since $S(\mathcal{A},\mathcal{A},\mathcal{B})>0$, there exist $s,t\in\mathcal{A}$ and $r\in\mathcal{B}$ such that 
    $ts=r$, Proposition \ref{prop:go}. Hence, $\langle r,s\rangle\simeq\mathbb{D}_{2n}$ and $\mathcal{A}$ contains 
    the subrack $Y=R\coprod S$ which is of type D by Lemma \ref{lem:D4n_tipoD}.
\end{proof}

\begin{exa}
\label{exa:Co1}
The conjugacy classes of involutions of the Conway group $Co_1$ are of type D.
In fact, $S(\textup{2A},\textup{2A},\textup{6E})=6$,
$S(\textup{2B},\textup{2B},\textup{6A})=2592$ and
$S(\textup{2C},\textup{2C},\textup{6A})=25920$.
\end{exa}

\begin{exa}
\label{exa:B}
The conjugacy class \textup{2C} of the Baby Monster group $B$ is of type D,
since $S(\textup{2C},\textup{2C},\textup{6C})=82752$.
\end{exa}

\subsection{The \textsf{ATLAS}}\label{subsect:from-atlas}

The \textsf{ATLAS} of Finite Groups, often simply known as the
\atlas, is a group theory book by John Conway, Robert Curtis,
Simon Norton, Richard Parker and Robert Wilson (with computational
assistance from J. G. Thackray), published in 1985 -- see
\cite{atlas}. It lists basic information about finite simple
groups such as presentations, conjugacy classes of maximal
subgroups, character tables and power maps on the conjugacy
classes.

The \textsf{ATLAS} is being continued in the form of an electronic
database -- see \cite{ATLASwww}. It currently contains information
(including $5215$ representations) on about $716$ groups. In order
to access to the information contained in the \atlas, we use the
\texttt{AtlasRep} package for \textsf{GAP} -- see
\cite{AtlasRep1.4}.

We recall some notations from the \atlas\/ for the reader not used
to it.

\noindent The notation for families of simple groups can be found in pages
x-xiii of \cite{atlas}:
\begin{itemize}
    \item
        $L_n(q) = \mathbf{PSL}(n,\mathbb{F}_q)$ is the projective special linear group.

        \smallbreak     \item $U_n(q) = \mathbf{PSU}(n,\mathbb{F}_q)$ is the projective
special unitary group.

\smallbreak     \item $S_{2n}(q) = \mathbf{PSp}(2n,\mathbb{F}_q)$ is the
projective symplectic group.

\smallbreak     \item $O_{n}(q) < \mathbf{PSO}(n,\mathbb{F}_q)$, for $n$ odd, and
$O_{2k}^{\epsilon}(q) < \mathbf{PSO}(2k,\mathbb{F}_q)$, for $n=2k$ even, is the
usual simple subgroup of the projective special orthogonal
group, where $\epsilon =\pm$ means the plus/minus type of the
corresponding quadratic form.

\smallbreak     \item $G_2(q)$, $E_6(q)$ are exceptional
groups in the family of Chevalley groups.

\end{itemize}

\bigbreak  There are various ways to combine groups or abbreviate
some groups structures -- see page xx of \cite{atlas}.
Assume that $K$ and $G$ are groups. Then:
\begin{itemize}
\item $K\ldotp
G$ means a group $L$ fitting into an extension $ 1\to  K \to L
\to  G \to   1$; which extension should be clear from the
context. Besides, $K : G$ means that the extension is split
(i.~e. when $K.G$ is a semi-direct product) and $K\cdot G$
means that the extension is non-split.

\item $K\times G$ is the direct product of $K$ and $G$.
\item $K^m$ denotes the direct product of $m$ groups isomorphic to
$K$.
\item $p^m$, for $p$ prime, denotes the elementary abelian group of order
$p^m$.
\item $[m]$, for $m\in\N$, denotes an abitrary group of order $m$.
\item $m$ denotes the cyclic group of $m$ elements.
\item $p^{n+m}$ indicates a case of $p^n.p^m$.
\item $p^{1+2n}$ or $p^{1+2n}_+$ or $p^{1+2n}_-$ is used for the particular case
    of an extraspecial group.
\end{itemize}

Product of three or more groups are left-associated. That is, $A.B.C$
means $(A.B).C$, and implies the existence of a normal subgroup
isomorphic to $A$.

\bigbreak We extracted from the \textsf{ATLAS}:
\begin{itemize}
\smallbreak
\item The relevant information about maximal subgroups of a given sporadic
    simple group $G$ and, if possible, the informa\-tion about the fusion of
    conjugacy classes from maximal subgroups of $G$ into $G$.

\item The notation for the conjugacy classes. The conjugacy classes that
    contain elements of order $n$ are named $n$A, $n$B, $n$C, \ldots, and
    notice that the alphabet used here is potentially infinite.
The conjugacy classes computed with \textsf{GAP} of a group given by a
particular representation (taken from the \textsf{ATLAS} or not) are not
necessarily named following the \textsf{ATLAS} notation. To avoid problems, in all these
cases, we print sizes of the centralizers to recognize the classes we
are working with.


\end{itemize}

\subsection{Notations and results from \textsf{GAP}}\label{subsect:from-gap}

We checked with \textsf{GAP} the following information:

\begin{itemize}
\bigbreak  \item
Real and quasi-real conjugacy classes. To recognize real conjugacy classes
we use the \textsf{GAP} function \texttt{RealClasses}. To recognize quasi-real conjugacy
classes we developed the \textsf{GAP} function \texttt{QuasiRealClasses}.
This function uses the \textsf{GAP} function \texttt{PowerMaps}. See \cite{logbook} for more details.

\bigbreak \item If $H$ is a subgroup of a group $G$, we use
the function \texttt{PossibleClassFusions} or
\texttt{FusionConjugacyClasses} for computing, or recovering,
the fusion of conjugacy classes from $H$ to $G$. See
\cite{logbook} for more details.

\bigbreak \item We computed with \textsf{GAP} sums
$S(\oc,\mathcal{K}_{1},\mathcal{K}_{2})$, where $\oc$ is the
conjugacy class of an involution and
$\mathcal{K}_{1},\mathcal{K}_{2}$ are conjugacy classes of
elements of order $3$, in order to apply Proposition
\ref{cor:ct_a4}. See \eqref{eqn:suma-a4}.

\end{itemize}

\section{On Theorem \ref{th:racks-liquidados}}\label{sect:proofs-b}

\subsection{Type D}\label{subsect:splitconj-rack}
Let $X$ be a rack, $n\in \N$, $\varGamma\subset\GL(n,\ku)$ a
subgroup,  $q:X\times X\to\varGamma$ a 2-cocycle. Let $V = \ku X\otimes\ku^{n}$ as in Subsection \ref{subsect:racks}
and $g: X \to \GL(V)$ be the morphism of racks given by
\begin{equation}\label{eqn:defgpo}
    g_x (e_yw) = e_{x\trid y}q_{x,y}(w), \qquad x,y\in X, w\in V.
\end{equation}

\begin{definition}\label{def:tiposbyd} Let $(X, \trid)$ be a rack and $q$ a 2-cocycle.

    \begin{itemize}
\item    We say that $X$ is \emph{faithful} if $\phi: X \to \s_X$ is
    injective.

\medbreak
\item We say that $(X, q)$ is
    \emph{faithful} if $q: X \to \GL(V)$ is injective;  if $X$ is clear from the
    context, we shall also say that $q$ is faithful. If $X$ is faithful, then $(X, q)$ is faithful for any $q$.

    \medbreak\item We say that $X$ \emph{collapses} if for any faithful cocycle $q$ with values
            in a finite group $\varGamma\subset\GL(n,\ku)$, for any $n\in \N$, $\dim\toba(X,q)=\infty$.

 \medbreak\item A rack $X$ is \emph{decomposable} iff there exist disjoint subracks
        $X_1,X_2\subset X$ such that
        $X=X_1\coprod X_2 $. Otherwise,  $X$ is \emph{indecomposable}.

\medbreak\item  We say that $X$ is \emph{of type D} if there
            exists  a decomposable subrack $Y = R\coprod S$ of
            $X$ such that
            \begin{equation}\label{eqn:hypothesis-subrack}
                r\trid(s\trid(r\trid s)) \neq s, \quad \text{for some } r\in R,
                s\in S.
            \end{equation}
    \end{itemize}
\end{definition}

Our main tool is the following rack-theoretical version of
\cite[Th. 8.6]{HS1}, whose proof uses the main result of
\cite{AHS}.

\begin{thm}\label{th:racks-claseD} \cite[3.6]{AFGV}.
If $X$ is a finite rack of type D, then $X$ collapses. \qed
\end{thm}

Let $X$ be a simple rack that is not a permutation rack. Then $X$ is
faithful, hence any cocycle is faithful, and the notion of
collapsing in Definition \ref{def:tiposbyd} and Definition
\ref{def-intro:collapses} in the Introduction coincide. Therefore,
for the purposes of this paper, we first need to check when a
conjugacy class in a sporadic simple group is of type D. We
collect some useful consequences of Theorem \ref{th:racks-claseD}.

\begin{lem}\label{lem:D_cocientes}
\renewcommand{\theenumi}{\roman{enumi}}\renewcommand{\labelenumi}{(\theenumi)}
\begin{asparaenum}

    \item Let $X$ be a rack of type D and let $Z$ be a finite rack that
        admits a rack epimorphism $f: Z\to X$, $X$ of type D.  Then $Z$
        is of type D.

\medbreak\item Let $f:G\to H$ be a group epimorphism and
let $g\in G$ such that $f(g)=h$. If the conjugacy class of $h$
is of type D, then the conjugacy class of $g$ is of type
D.
    \end{asparaenum}\end{lem}

\pf For (i), $\pi^{-1}(Y) = \pi^{-1}(R)\coprod \pi^{-1}(S)$ is a decomposable
        subrack of $Z$ satisfying \eqref{eqn:hypothesis-subrack}. Clearly, (ii) follows from (i).
\epf

Notice the inference of Theorem \ref{th:racks-liquidados-aplicacion}
from Theorem \ref{th:racks-liquidados} by Lemma \ref{lem:D_cocientes} (1).

\begin{thm}\label{th:racks-liquidados-aplicacion}
Let $G$ be any finite group,  $\Q$ a conjugacy class of $G$ and
$g\in \Q$. Assume that there is a rack epimorphism $\Q\to \oc$,
where $\oc$ is a  non-trivial conjugacy class of a sporadic group
NOT listed in Table \ref{tab:0}. Then $\dim \toba(\Q, \rho) =
\infty$ for every $\rho\in \Irr C_{G}(g)$. \qed
\end{thm}

\smallbreak
Let us say that a finite group is \emph{of type D} if all its non-trivial conjugacy classes
are of type D. By Th. \ref{th:racks-claseD}, a finite group of type D collapses.

\begin{prop}
Let $0\to K\overset{i}{\to}G\overset{p}{\to}H\to0$
be a short exact sequence of finite groups such that $K$ and $H$ are of
type D. Then $G$ is of type D. \end{prop}
\begin{proof}
Let $g\in G$. If $p(g)=1$, then $g\in i(K)$ and the conjugacy class
of $g$ in $G$ is of type D (since the conjugacy class in $K$ of
$k\in K$, where $i(k)=g$, is of type D). If $1\ne p(g)$, then the
conjugacy class of $h=p(g)$ in $H$ is of type D and then, by Lemma
\ref{lem:D_cocientes}, the conjugacy class of $g$ in $G$ is of
type D.
\end{proof}

We get immediately the following result.

\begin{thm}\label{th:factores-colapsan}
Let $G$ be any finite group. Assume that the simple factors of $G$ in the Jordan-H\"older decomposition
of $G$ are of type D. Then $G$ is of type D, hence it collapses. \qed
\end{thm}

By Theorem \ref{th:racks-liquidados}, the groups $Th$, $He$ and $HN$ are of
type D; we also know that the groups $G_2(3)$, $G_2(5)$ are of type D.

\bigbreak
\subsection{Algorithms}\label{subsect:algorithms}

We now explain our algorithms to implement the technique of the previous subsection.

\begin{alg}\label{algo:racks-claseD} Let $\varGamma$ be a
finite group and let $\oc$ be a conjugacy class. Fix $r\in \oc$.

\begin{enumerate}
    \item For any $s\in \oc$, check if $(rs)^2 \neq (sr)^2$; this
    is  equivalent to
    \eqref{eqn:hypothesis-subrack}.

    \smallbreak
    \item If such $s$ is found, consider the subgroup $H$ generated
    by $r$, $s$. If $\oc^H_r \cap \oc^H_s =
    \emptyset$, then $Y = \oc^H_r \coprod \oc^H_s$ is
    the decomposable subrack we are looking for and $\oc$ is of
    type D.
\end{enumerate}
\end{alg}

We have found that a useful variant, instead of going over all the elements, is to
choose randomly a certain number of $s\in \oc$ and check the conditions above.
This turns out to be very often a much quicker way to see if $\oc$ is of type D.

\medbreak In the practice, for large groups, it is more economical to
implement the algorithm in a recursive way.

\begin{alg}\label{algo:racks-claseD-conmaximales}
Let $G$ be a finite group.

\begin{enumerate}
    \item List all maximal subgroups of $G$ up to conjugacy, say
    $\M_1, \dots, \M_k$, with $\vert \M_1\vert \leq \vert \M_2\vert \dots \leq \vert \M_k\vert$.

    \smallbreak
    \item Perform the Algorithm \ref{algo:racks-claseD} for every
    conjugacy class of $\varGamma = \M_1$. Let $\D_1$ be the set of
    conjugacy classes of $G$ that contain a conjugacy class of $\M_1$ of type D.

\smallbreak
    \item Perform the Algorithm \ref{algo:racks-claseD} for every
    conjugacy class of $\varGamma = \M_2$ that is not contained in
    any $\oc \in \D_1$. Let $\D_2$ be the set of
    conjugacy classes of $G$ that either contain a conjugacy class of $\M_2$ of type D, or else
    are in $\D_1$.

\smallbreak
    \item Continue in this way, performing the Algorithm \ref{algo:racks-claseD} for every
    conjugacy class of the various maximal subgroups $\M_j$ and producing at each step a set of discarded
    classes $\D_j$.

\smallbreak
    \item Perform the Algorithm \ref{algo:racks-claseD} for every
    conjugacy class $\oc$ of $G$ not in $\D_k$. Let $\D$ be
    the set of conjugacy classes of $G$ that either are of
    type D by this argument, or else are in $\D_k$. This set
    $\D$ is the output of the algorithm.
\end{enumerate}
\end{alg}

The treatment of a maximal subgroup $\M$ can be simplified if it fits into a
short exact sequence of groups $1\to K\overset{i}{\to}\M\overset{p}{\to}H\to 1$,
where we know the conjugacy classes in $H$ of type D. To apply Lemma \ref{lem:D_cocientes} to
a  conjugacy class in $\M$, we just need to know that some specific conjugacy classes in $H$ are of type D.

\begin{lem}
\label{lem:criterioD}
Let $g\in \M$ of order $m$.
Assume that every conjugacy class in $H$ with representative of order $k$ is of type $D$,
for every $k$ such that $k\mid m$ and $\frac{m}{k}\mid \vert K\vert$.
Then the conjugacy class of $g$ in $\M$ is also of type D.
\end{lem}

\begin{proof}
If  $h=p(g)$ has order $k$, then $k\mid m$.
Also, since $g^{k}\in\ker(p)=i(K)$ and $i$ is a monomorphism, we
have that $\frac{m}{k}=\frac{m}{\left(m,k\right)}\mid \vert K\vert$. Thus, it suffices
to have these conjugacy classes of type D to
conclude from Lemma \ref{lem:D_cocientes}. \end{proof}

Here is another useful remark.

\begin{lem}
\label{lem:Ddirecto}
Let $G=H\times K$ be a direct product of finite groups. Let $\mathcal{O}$
be a conjugacy class of $H$ with representative $h$, and $k\in K$.
If $\mathcal{O}$ is of type D, then the conjugacy class of $h\times k$
is of type D.\qed
\end{lem}

The proof of Theorem \ref{th:racks-liquidados} follows by applying Algorithm
\ref{algo:racks-claseD-conmaximales} either to the sporadic groups or their
maximal subgroups. The details of the computations are in \cite{logbook},
except for the involutions treated in Examples \ref{exa:Co1} and \ref{exa:B}.

\section{On Theorem \ref{teor:complete}}
\label{sect:proofs-a}

\subsection{Proof of Theorem \ref{teor:complete}}\label{subsect:proofs-a}

\begin{table}[ht]
\caption{Proof of Theorem \ref{teor:complete}.}\label{tab:final}
\begin{center}
\begin{tabular}{|p{1cm}|c|c|}
\hline $G$ & {\bf Classes}  & {\bf Relevant information}
\\ \hline  $M_{11}$ & \textup{8A, 8B} & Example \ref{exa:affin-notipoD}
\\ \cline{2-3} & \textup{11A, 11B} & \small{quasi-real} $j=3$, $g^{9}\neq g$
\\ \hline  $M_{12}$ &  \textup{11A, 11B} &  \small{quasi-real} $j=3$, $g^{9}\neq g$
\\ \hline  $M_{22}$ &  \textup{11A, 11B} &  \small{quasi-real} $j=3$, $g^{9}\neq g$
\\ \hline  $M_{23}$ &  \textup{23A, 23B} &  \small{quasi-real} $j=2$, $g^{4}\neq g$
\\ \hline  $M_{24}$ &  \textup{23A, 23B} &  \small{quasi-real} $j=2$, $g^{4}\neq g$
\\ \hline  $J_{2}$ &  \textup{2A} &  $S(2A,3B,3B)=18$
\\ \cline{2-3} &  \textup{3A} & real
\\ \hline  $Suz$ &  \textup{3A} &  real
\\ \hline  $HS$ &  \textup{11A, 11B} &   \small{quasi-real} $j=3$, $g^{9}\neq g$
\\ \hline  $McL$ &  \textup{11A, 11B} &   \small{quasi-real} $j=3$, $g^{9}\neq g$
\\ \hline  $Co_{3}$ &  \textup{23A, 23B} & \small{quasi-real} $j=2$, $g^{4}\neq g$
\\ \hline  $Co_{2}$ &  \textup{2A} & abelian subrack
\\ \cline{2-3} & \textup{23A, 23B} &  \small{quasi-real} $j=2$, $g^{4}\neq g$
\\ \hline  $Co_{1}$ &   \textup{3A} &   real
\\ \cline{2-3} &  \textup{23A, 23B} &   \small{quasi-real} $j=2$, $g^{4}\neq g$
\\ \hline  $J_{1}$ &  \textup{15A, 15B, 19A, 19B, 19C} & real
\\ \hline  $O'N$ &  \textup{31A, 31B}  & \small{quasi-real} $j=2$, $g^{4}\neq g$
\\ \hline $J_{3}$ &  \textup{5A, 5B} & real
\\ \cline{2-3} &  \textup{19A, 19B} &   \small{quasi-real} $j=4$, $g^{16}\neq g$
\\ \hline  $Ly$    &  \textup{33A, 33B} &   \small{quasi-real} $j=4$, $g^{16}\neq g$
\\ \cline{2-3} &  \textup{37A, 37B, 67A, 67B, 67C} &   real
\\ \hline  $Ru$ &  \textup{29A, 29B} & real
\\ \hline   $J_4$ &  \textup{ 29A, 37A, 37B, 37C, 43A, 43B, 43C} &  real

\\ \hline  $Fi_{22}$   & \textup{2A} & abelian subrack
\\ \hline    $Fi_{23}$ & \textup{2A} & abelian subrack
\\ \cline{2-3} & \textup{23A, 23B} &   \small{quasi-real} $j=2$, $g^{4}\neq g$

\\ \hline  $Fi'_{24}$ &  \textup{27B, 27C, 29A, 29B} & real

\\ \cline{2-3} &   \textup{33A, 33B, 39C, 39D} & real

\\ \cline{2-3} &   \textup{23A, 23B} & \small{quasi-real} $j=2$, $g^{4}\neq g$

\\ \hline   $B$ &  \textup{2A} & abelian subrack

\\ \cline{2-3} &
\textup{47A, 47B} & \small{quasi-real} $j=2$, $g^{4}\neq g$
\\ \hline   $M$ &  \textup{41A}& real

\\ \cline{2-3} &   \textup{47A, 47B, 69A, 69B, 71A, 71B, 87A, 87B} & \small{quasi-real} $j=2$, $g^{4}\neq g$

\\ \cline{2-3} &   \textup{59A, 59B} & \small{quasi-real} $j=3$, $g^{9}\neq g$

\\ \hline $T$ &  \textup{2A} & $S(2A,3A,3A)=108$
\\ \hline
\end{tabular}
\end{center}
\end{table}

As explained in Subsect. \ref{subsect:yd}, Theorem \ref{teor:complete} is
equivalent to the following statement.

\begin{thm}\label{teor:complete-nichols}
    \
\begin{enumerate}
\item
If $G$ is a sporadic group but not $Fi_{22}$, $B$, $M$, then the Nichols algebra of any  irreducible Yetter-Drinfeld
module has infinite dimension.


\item If $G$ is $Fi_{22}$, $B$ or $M$ and the pair $(\oc,\rho)$ is not listed in Table 1, then
$\dim\toba(\mathcal{O},\rho)=\infty$.
\end{enumerate}
\end{thm}

\pf By Theorems \ref{th:racks-liquidados} and \ref{th:racks-claseD}, it remains to consider
the conjugacy classes listed in Table \ref{tab:0}; we summarize the methods for each class in Table \ref{tab:final}.
The corresponding Nichols algebras have infinite dimension by the following reasons:

\begin{itemize}
    \item If the class is real, then Lemma \ref{lem:inversos} applies.

\medbreak\item If the class is quasi-real, then Lemma  \ref{lem:potencias} applies.

\medbreak\item If the class contains an involution, then often Proposition
\ref{cor:ct_a4} applies, since some sum (shown in the Table) is not zero.

\medbreak\item The conjugacy classes labeled \textup{2A} of
the groups $Co_2$, $Fi_{22}$, $Fi_{23}$ or $B$ do not
collapse, but $\dim\mathfrak{B}(\mathcal{O}, \rho)=\infty$ for
any irreducible representation $\rho$ of the corresponding
centralizer. Indeed, it is enough to consider the
representations $\rho$ such that $\rho(g)=-1$. By Lemma
\ref{lem:contar}, we are reduced to find an involution $x$
such that $gh=hg$, for $h=xgx^{-1}$, and to compute the
multiplicities of the eigenvalues of $\rho(h)$. For this last
task, we use Remark \ref{rem:contar_autovalores}. In the case
of the group $Fi_{23}$,  we first apply Lemma \ref{lem:contar}
to $S_8(2)$, that is isomorphic to a subgroup of $Fi_{23}$,
and then apply Lemma \ref{lema:technique-fusionsubgroups}.
Finally, the class labeled \textup{2A} of $Fi_{23}$ embeds as
a subrack of the class labeled \textup{2A} of $B$, hence the
technique of the subgroup applies. See \cite{logbook} for
details.
\end{itemize}

We finally explain the representations appearing in
Table \ref{tab:1}.

Let $\oc$ be one of the conjugacy classes \textup{22A, 22B} of
$Fi_{22}$, \textup{46A, 46B} of $B$,  \textup{92A, 92B, 94A, 94B}
of $M$. Then the centralizer of $g\in \oc$ is the cyclic group
$\langle g\rangle$. Thus $\Irr C_{G}(g)$ is also cyclic; say
$\chi_\omega$ is the representation such that $\chi_\omega(g) =
\omega$, $\omega$ a root of 1 whose order divides the order of
$g$. On the other hand these 
classes are quasi-real, as stated in Table \ref{tab:1}. By Lemmata
\ref{lem:inversos} or \ref{lem:potencias}, only $\chi_{-1}$
survives.

Let $\oc$ be the conjugacy class \textup{34A} of $B$ and $g \in \oc$. This
class is real. By Lemma \ref{lem:inversos}, we discard all irreducible
representations $\rho$ of the centralizer $\Z/34 \times \Z/2$ except those
satisfying $\rho(g)= -1$, i.~e.  $\rho=\chi_{-1}\otimes \epsilon$,
$\chi_{-1}\otimes \sgn$, where $\epsilon$ and $\sgn$ mean the trivial and the
sign representation of $\Z/2$.

Let $\oc$ be one of the conjugacy classes \textup{46A, 46B} of $M$  and $g \in
\oc$. The class $\oc$ is quasi-real, as stated in Table \ref{tab:1}. By Lemma
\ref{lem:potencias}, we discard all irreducible representations $\rho$ of the
centralizer $\Z/23 \times \mathbb D_4$ except the one satisfying $\rho(g)= -1$,
i.~e.  $\rho=\epsilon \otimes \rho_2$, where $\epsilon$ is the trivial
representation of $\Z/2$ and $\rho_2$ is the unique irreducible representation
of $\mathbb D_4$ of degree 2.
\epf


In the proof of Theorem \ref{teor:complete-nichols}, we needed the structure of the centralizers of some specific elements in some sporadic groups; this was kindly communicated to us by Thomas Breuer, when not available in the literature.

\subsection{Remarks on the remaining irreducible Yetter-Drinfeld modules in Table \ref{tab:1}}\label{subsect:remarks-remaining}

Let $G$ be a sporadic group and $\oc$ a conjugacy class as in  Table \ref{tab:1}.
Assume that $r$ and $s$ are two elements in $\oc$ such that \eqref{eqn:hypothesis-subrack} holds
and let $H$ be the subgroup of $G$ generated by $r$
and $s$. If the conjugacy classes in $H$ of $r$ and $s$ are disjoint, then $\oc$ would be
of type D. Notice that $H$ should be contained in a maximal subgroup $\M$ of $G$ and, clearly, it would be enough to perform the
necessary computations in $\M$. Also,
if $G$ is not the Monster group, the list of all maximal subgroups is
known. So, we proceed to investigate these maximal subgroups by the fusion of
conjugacy classes, see Subsection \ref{subsect:from-gap}.


\subsubsection{The classes $\oc$ labeled \textup{22A, 22B} in $Fi_{22}$}

We know that $H$ should be contained in $\M_1 \simeq 2.U_6(2)$. 
It is not possible with our computational resources to determine
if these classes are of type D in $\M_1$, i.~e. due to the
size of these conjugacy classes it is not possible to complete all
the possible elections of $r$ and $s$. Actually, in all the
computed cases $H$ yields to be the centralizer of the
representative of the class or $\M_1$ itself.

The same occurs if we want to determine if the conjugacy classes
\textup{11A, 11B} in $U_6(2)$ are of type D. So we cannot
decide that the corresponding Nichols algebra is
infinite-dimensional by lifting racks of type D from the classes \textup{11A, 11B} of $U_6(2)$.



\subsubsection{The classes $\oc$ labeled \textup{16C, 16D, 32A, 32B, 32C, 32D, 34A} in $B$}

The classes \textup{16C, 16D} meet the following maximal
subgroups:
\begin{align*}
\M_1 &\simeq 2.\hspace{2pt}{}^2\hspace{-2pt}E_6(2):2, & \M_2 &\simeq 2^{1+22}.Co_2,\\
\M_4 &\simeq 2^{9+16}.S_8(2), &\M_7 &\simeq 2^{2+10+20}.(M_{22}:2 \times \s_3),\\
\M_8 &\simeq [2^{30}].L_5(2), & \M_{10} &\simeq [2^{35}].(\s_5\times L_3(2));
\end{align*}
the classes \textup{32A, 32B} meet $\M_2$, $\M_4$, $\M_7$ or
$\M_{10}$; the classes \textup{32C, 32D} meet $\M_1$, $\M_2$,
$\M_4$, $\M_7$ or $\M_{10}$; the class \textup{34A} meets $\M_1$.
The known matrix representations of these subgroups do not allow
us to perform the necessary computations.  Notice that the sixth
maximal subgroup of $B$ is excluded of our analysis, because we do
not know the fusion of conjugacy classes $\M_6\to B$.


\subsubsection{The classes $\oc$ labeled \textup{46A, 46B} in $B$}

We know that $H$ should be contained in the second maximal subgroup $\M_2$.  As
already said, the computations for this case are out of our resources.  Note
also that $Co_2$ has no elements of order 46; thus, if $g\in \M_2$ has order
46, then its projection in $Co_2$ has order 23, but the classes of elements of
order 23 in $Co_2$ are not known to be of type D.

\subsubsection{The classes $\oc$ labeled \textup{32A, 32B, 46A, 46B, 
92A, 92B, 94A, 94B} in $M$}

Among the known maximal subgroups of the Monster, \cite{B} provides the fusion of only 4 of them:
$\M_1 \simeq 2.B$, $\M_2 \simeq 2^{1+24}.Co_1$, $\M_3 \simeq 3.Fi_{24}$, and $\M_{9} \simeq \s_3 \times Th$.
The classes \textup{32A, 32B} meet $\M_{1}$, $\M_{2}$; the classes \textup{46A, 46B} meet $\M_{1}$, $\M_{2}$;
the classes \textup{92A, 92B} meet $\M_{2}$; the classes \textup{94A, 94B} meet $\M_{1}$. As before,
the necessary computations are out of our resources.


\begin{acknowledgements*}
Some of the results of the present paper were announced in
\cite{AFGV3} and at several conferences like: 
IV Encuentro Nacional
de \'Algebra, August 2008, La Falda, Argentina;
Grou\-pes quantiques dynamiques et cat\'egories de
fusion, April 2008, CIRM, Luminy, France; 
First De Br\'un Workshop on Computational Algebra,
August 2008, Galway, Ireland; Hopf algebras and related topics (A conference in honor of Professor
Susan Montgomery), February 2009, University of Southern
California; XVIII Latin American Algebra Colloquium, August 2009, S\~ao Pedro, Brazil; Coloquio de
\'Algebras de Hopf, Grupos Cu\'anticos y Categor\'{\i}as Tensoriales,
August, 2009, La Falda, Argentina; XIX Encuentro
Rioplatense de \'Algebra y Geometr\'{\i}a Algebraica, November
2009, Montevideo, Uruguay; VII workshop in Lie Theory and its applications, November
2009, C\'ordoba, Argentina. We are grateful to the organizers for the kind invitations.

We thank Sebasti\'an Freyre for interesting conversations. We are
very grateful to Alexander Hulpke, John Bray, Robert Wilson  and very specially to
Thomas Breuer for answering our endless questions on \textsf{GAP}.
We also thank Enrique Tobis and the people from \textsf{shiva} and \textsf{ganesh}
who allowed us to use their computers.
\end{acknowledgements*}


\begin{thebibliography}{AFGV2}


\bibitem[A]{A}{\sc N. Andruskiewitsch},
\emph{Some remarks on Nichols algebras}.
In "Hopf algebras", Bergen, Catoiu and Chin (eds.), 25-45 (2004), M. Dekker.

\bibitem[AF]{AF2} {\sc N. Andruskiewitsch and F. Fantino},
\emph{On pointed Hopf algebras associated with alternating and
dihedral groups}, Rev. Uni\'on Mat. Argent. 48-3, (2007), 57-71.


\bibitem[AFGV1]{AFGV} {\sc N. Andruskiewitsch,} {\sc F. Fantino}, {\sc M.
Gra\~na} {\sc and  L. Vendramin}, \emph{Finite-di\-mensional
pointed Hopf algebras with alternating groups are trivial}.
Ann. Mat. Pura Appl, \verb+doi:10.1007/s10231-010-0147-0+.

\bibitem[AFGV2]{logbook}\bysame,
\emph{The logbook of pointed Hopf algebras over the sporadic simple groups}, preprint:\texttt{arXiv:1001.1113}.

\bibitem[AFGV3]{AFGV3} \bysame, \emph{Pointed Hopf algebras over
some sporadic simple groups}, 
C. R. Math. Acad. Sci. Paris 348 (2010) pp. 605-608


\bibitem[AG]{AG1} {\sc N. Andruskiewitsch and  M. Gra\~na},
\emph{From racks to pointed Hopf algebras}, Adv. Math.
\textbf{178}  (2003), 177 -- 243.

\bibitem[AHS]{AHS} {\sc N. Andruskiewitsch, I. Heckenberger and H.-J. Schneider},
{\em The Nichols algebra of a semisimple Yetter-Drinfeld module},
\texttt{arXiv:0803.2430}.

\bibitem[AS1]{AS-cambr}  {\sc N. Andruskiewitsch and  H.-J. Schneider},
\emph{Pointed Hopf Algebras}, in ``New directions in Hopf
algebras'', 1--68, Math. Sci. Res. Inst. Publ. \textbf{43},
Cambridge Univ. Press,  2002.

\bibitem[AS2]{AS-ann}\bysame, \textit{On the
classification of finite-dimensional pointed Hopf algebras}. 
Ann. Math. Vol. 171 (2010), No. 1, 375417.


\bibitem[AZ]{AZ}  {\sc N. Andruskiewitsch and  S. Zhang}, \emph{On pointed Hopf algebras associated to some conjugacy
classes in $\mathbb S_n$}, Proc. Amer. Math. Soc. 135  (2007),
2723-2731.

\bibitem[B]{B}
{T. Breuer}, \emph{The GAP Character Table Library, Version 1.2 (unpublished)}; 
\texttt{http://www.math.rwth-aachen.de/\~{}Thomas.Breuer/ctbllib/}.

\bibitem[CC+]{atlas} {\sc J. H. Conway, R. T. Curtis, S. P. Norton, R. A. Parker and R. A.
Wilson}, \emph{Atlas of finite groups}, Oxford University Press,
1985.

\bibitem[D]{dickson} {\sc  L. Dickson}, \emph{The alternating group on eight letters and the quaternary
              linear congruence group modulo two},
Math. Ann. {\bf 54} (1901), pp. 564--569.

\bibitem[F]{fantino-2007} {\sc  F. Fantino}, \emph{On pointed Hopf algebras associated with
Mathieu groups}, J. Algebra Appl. 8 5 (2009) 633-672.


\bibitem[FGV]{FGV}  {\sc S. Freyre,  M. Gra\~na and  L. Vendramin},
\emph{On Nichols algebras over $\mathbf{GL}(2,\mathbb{F}_q)$ and
${\mathbf{SL}(2,\mathbb{F}_q)}$}, J. Math. Phys. \textbf{48}
(2007), 123513-1 -- 123513-11.


\bibitem[GAP]{GAP}
 The \textsf{GAP}~Group, \emph{\textsf{GAP} -- Groups, Algorithms, and Programming,
 Version 4.4.12};
 2008,
Available at \verb+http://www.gap-system.org+.

\bibitem[Go]{g}  {\sc D. Gorenstein}, \emph{Finite groups},
Harper \& Row Publ., New York, (1968), xv+527.



\bibitem[H]{H-all} {\sc I. Heckenberger},
\textit{Classification of arithmetic root systems},
Adv. Math. {\bf 220} (2009) 59--124.


\bibitem[HS]{HS1} {\sc I. Heckenberger and  H.-J. Schneider},
\emph{Root systems and Weyl groupoids for semisimple Nichols
algebras}. 
Proc. London Math. Soc. \verb+doi:10.1112/plms/pdq001+.




\bibitem[L]{Lu} {\sc G. Lusztig}, \emph{Introduction to quantum groups}, Birkh\"auser,
1993.

\bibitem[R]{Ro} {\sc M. Rosso},
\textit{Quantum groups and quantum shuffles}, Invent. Math.
\textbf{133 } (1998),  399--416.

\bibitem[S]{Sbg} {\sc P.~Schauenburg}, {\it A characterization of the Borel-like
subalgebras of quantum enveloping algebras}, Comm. in Algebra
{\bf24} (1996), pp. 2811--2823.

\bibitem[WWT+]{ATLASwww} {\sc R. A. Wilson, P. Walsh, J. Tripp, I. Suleiman,
R. Parker, S. Norton, S. Nickerson, S. Linton, J. Bray and R. Abbott}, \emph{A world-wide-web Atlas of finite group representations},
\texttt{http://brauer.maths.qmul.ac.uk/Atlas/v3/}.

\bibitem[WPN+]{AtlasRep1.4} {\sc R. A. Wilson, R. A. Parker, S. Nickerson,
J. N. Bray and T. Breuer}, \emph{{AtlasRep}, A \textsf{GAP} Interface to the \textsf{ATLAS} of
Group Representations,
{V}ersion 1.4}, 2007, Refereed \textsf{GAP} package, 
\verb+http://www.math.rwth-aachen.de/~Thomas.Breuer/atlasrep+.

\end{thebibliography}
\end{document}